\documentclass[a4paper,11pt]{article}

\usepackage{amsmath}
\usepackage{booktabs}
\usepackage{graphicx}
\usepackage[english]{babel}
\usepackage{nicefrac}
\usepackage{subcaption}

\title{Na\"{i}ve Newsvendor Adjustments:\\
Are They Always Detrimental?}

\author{Congzheng Liu\thanks{Department of Management Science,
Lancaster University, Lancaster LA1 4YX, UK.
Email: {\tt \{c.liu19,a.n.letchford,i.svetunkov\}@lancaster.ac.uk}}
\and Adam N.\ Letchford$^*$ \and Ivan Svetunkov$^*$} 

\date{July 2022}

\begin{document}
\maketitle

\begin{abstract}
Newsvendor problems are an important and much-studied topic in stochastic
inventory control. One strand of the literature on newsvendor problems is
concerned with the fact that practitioners often make judgemental
adjustments to the theoretically ``optimal" order quantities. Although
judgemental adjustment is sometimes beneficial, two specific kinds of
adjustment are normally considered to be particularly na\"{\i}ve:
\emph{demand chasing} and \emph{pull-to-centre}. We discuss how these
adjustments work in practice and what they imply in a variety of settings.
We argue that even such na\"{\i}ve adjustments can be useful under certain
conditions. This is confirmed by experiments on simulated data. Finally,
we propose a heuristic algorithm for ``tuning” the adjustment parameters
in practice.
\\*[2mm]
{\bf Keywords:} Inventory; Judgemental adjustments; Demand chasing;
Pull-to-centre effect
\end{abstract}


\section{Introduction}

Single-period stochastic inventory control problems, known as
\emph{Newsvendor problems} (NVPs), have received much attention from the
Operational Research community (see, e.g., the books
\cite{Ch12,HW63,SPP98,Zi00}). In this paper, we focus on the simplest NVP,
in which there is only one product. The demand for the product over the
selling period is a random variable $\tilde d$, with known distribution. The
product is purchased before the period at a fixed unit price $v$, and sold
during the period at a unit price $p$. If any excess stock remains at the end
of the period, a holding cost $c_h$ is incurred per unit. If there is any
unsatisfied demand, a shortage cost $c_s$ is incurred per unit.

For a given $x$ and a given realisation $d$ of $\tilde{d}$, the realised profit
$\pi$ over period is:
\begin{equation} 
    \pi(x,d)=
    \begin{cases}
        pd-vx-c_h(x-d),& \text{if } x\geq d\\
        px-vx-c_s(d-x),& \text{if } x< d.
    \end{cases}
\end{equation}
The optimal order quantity that maximises the expected profit is then (\cite{AHM51,Ch12}):
\begin{equation} \label{eq:textbook}
    x^* = F^{-1}( \tau),
\end{equation}
where $F$ is the cumulative distribution function for $\tilde d$,
$c_o = v+c_h$ is the \emph{overage} cost, $c_u = p-v+c_s$ is the
\emph{underage} cost, and the quantile $\tau=\nicefrac{c_u}{(c_o+c_u)}$.
We will call a product with $\tau > 0.5$ a ``high-margin" product, and a
product with $\tau < 0.5$ a ``low-margin" product.

In the textbook formulation, it is assumed that one has a correct model of
the demand distribution $F$, with correct parameters. In real life, however,
correctness is rarely assured and a model typically misses important
information. Moreover, even if the model is correct, the parameters may
evolve over time (for example, due to market shocks or product innovations
by competitors). For these reasons, decision makers often make so-called
``judgemental adjustments" to the theoretically ``optimal" order quantities.
The literature on this topic is extensive (e.g.,
\cite{Be08,BK08,BOT12,BHS08,LB13,LHB14,SC00}).

Although judgemental adjustments can exist in many forms (see
\cite{OFK17,GW14,Ka14} for example), two specific kinds of adjustment have
received most attention in the NVP literature
(\cite{Be08,BK08,BHS08,LB13,LHB14}). The first, called \emph{pull-to-centre},
means adjusting the order quantity towards the mean demand. The second,
called \emph{demand chasing}, means adjusting it towards the demand in the
immediately preceding period. Although these two kinds of adjustment are
regarded as especially na\"{\i}ve, considerable evidence for their existence
has been found by scholars, and several theories have been developed to
explain them (\cite{Be08,BK08,BHS08,FKZ11,MHD13,WN14}). Yet, to our
knowledge, there has not been any numerical study of the effect that these
two na\"{\i}ve adjustment mechanisms are likely to have on the long-term
expected profit in the NVP.

In this paper, we look at na\"{\i}ve adjustment from a different point of
view. We begin by arguing that such adjustments may be useful under certain
conditions, due to the fact that any given statistical model of the demand
is unlikely to be completely accurate. In particular, we suggest that a
modest amount of na\"{\i}ve adjustment may be beneficial in two specific
situations that are highly important in practice; namely, (i) when a
relatively small amount of demand data is available, and (ii) when the
true demand model is unknown. This idea is tested via extensive
computational  experiments. We then consider the possibility of applying
na\"{\i}ve adjustments in an automated fashion. For this purpose, we
propose a simple heuristic for tuning the adjustment parameters. Finally,
we test the heuristic on a real-life example, with encouraging results.

The paper is organised as follows. In Section \ref{se:lit}, we summarise
the existing literature. In Section \ref{se:alt}, we argue in favour of
na\"{\i}ve adjustment, and propose a mathematical model which incorporates
both types of adjustment. In Section \ref{se:exp}, we conduct
experiments on simulated data under different conditions, and discuss the
results. In Section \ref{se:algorithm}, we present the tuning heuristic
and apply it to a real-life example. Finally, Section \ref{se:end} contains
some concluding remarks.


\section{Literature Review} \label{se:lit}

In this section, we review and discuss the existing literature on na\"{\i}ve newsvendor
adjustment. Subsections \ref{sub:dc} and \ref{sub:ptc} are concerned with
demand-chasing and pull-to-centre, respectively. 

\subsection{Demand chasing effect} \label{sub:dc}

The \emph{demand chasing} effect (DC), in the NVP context, describes a
phenomenon whereby decision makers tend to adjust their order quantity
towards the realised demand in the previous operating period. 

To our knowledge, the first paper to give empirical evidence for DC was
\cite{SC00}. They observed decision makers over fifteen consecutive
ordering periods, for a selection of products, each with known distribution.
They showed that the participants systematically deviated from the optimal
order quantity. They also gave a tentative explanation for the phenomenon,
based on the avoidance of \emph{regret}. The idea is that, if decision
makers fail to choose the \emph{ex-post} optimal order quantity in a given
period, they regret their decision, which leads them to adjust the
quantity in the next period.

After the publication of \cite{SC00}, evidence for DC appeared in many
papers. \cite{BK08} repeated the experiment, but with 100 decision periods.
They found that the participants tended to improve over time, but only very
slowly. \cite{Be08} studied the convergence of the participants behaviour,
and argued that the order quantities from decision makers converge to a
level different from the one which optimises expected profit. Other relevant
works include \cite{BHS08,Cu13,FKZ11,MHD13,WN14}.

\cite{LB13} argued
that some of the statistical techniques used in the behavioural experiments
were flawed. See the recent paper \cite{KM20} for a discussion of this issue.

Several models of DC were proposed in the above-mentioned papers. For brevity,
we present only the simplest model, which appeared in \cite{BHS08}. It takes
the form:
\begin{equation} \label{eq:dc}
    x_t = x_{t-1} + \beta (d_{t-1}-x_{t-1}),
\end{equation}
where $x_t$ is the actual order quantity in time period $t$, $d_{t-1}$ is the
realised demand in the previous period, and $\beta > 0$ is the
\emph{DC parameter}. A higher $\beta$ indicates a stronger demand chasing
effect. (In rare cases, one might observe $\beta < 0$, as a ``pull forward
in demand".)

We remark that the model \eqref{eq:dc} is equivalent to ``simple
exponential smoothing" or SES (\cite{Br56}), in the so-called
``error-correction" form. We observe that it ignores the NVP solution and
asymptotically converges to the mean demand. Moreover, it assumes a fixed
$\beta$ for all periods, implying that practitioners adjust the orders
every period by the same proportion. Fortunately, the latter assumption
does not cause serious problems in practice. Indeed, even if practitioners
adjust the order with different quantities over time, their behaviour can
be modelled on average using \eqref{eq:dc}. Furthermore, there is empirical
evidence that, for human decision makers, the value of their $\beta$ will
eventually converge to a single value over time (\cite{SC00,ZS19}).


\subsection{Pull-to-centre effect} \label{sub:ptc}

The \emph{pull-to-centre} effect (PtC), also known as the mean anchor
heuristic, describes a phenomenon when a decision maker adjusts the order
quantity towards the mean demand.

To our knowledge, the first paper to give empirical evidence for PtC was
again \cite{SC00}. Their interpretation of PtC is that decision makers order
less than the optimal amount for high-margin products, but more for
low-margin products. They also discussed several possible causes for the
phenomenon, including risk and loss aversion, underestimation of opportunity
cost, and waste aversion. They also discussed a possible explanation in
terms of ``prospect theory" (\cite{KT79}).

Alternative explanations of PtC include adaptive learning (\cite{Be08}),
decision noise and optimisation error (\cite{Su08}), overconfidence bias
(\cite{RC13}), and psychological costs associated with leftovers and
stockouts (\cite{HLC10}). We remark that \cite{LHB14} argued that some of
the statistical techniques used to detect PtC were flawed, just as
\cite{LB13} argued for DC.

There is some evidence that individual differences can affect the behaviour
of the decision maker in NVPs. \cite{De13} showed that males tend to take
more risks than females, which leads them to order more, on average.
\cite{Cu13} and \cite{FKZ11} showed that differences in nationality
correlate with different biases while making newsvendor decisions.

\cite{Be08} found that decision makers tend to be more biased towards the
mean demand in earlier periods than in later periods. This suggests that
training could be of some benefit to the decision making. Additional
discussions of training effects can be found in \cite{BK08,BHS08,RC13,ZS19}.

Following the works of \cite{Be08} and \cite{BHS08}, the PtC effect can be
expressed mathematically as:
\begin{equation} \label{eq:PtC}
    x_t = (1-\gamma) x_t^* + \gamma \hat{\mu}_t,
\end{equation}
where $\hat{\mu}_t$ is the estimated mean demand for the period $t$ (which
can be obtained with a forecasting technique), and $0<\gamma<1$ is the
\emph{PtC parameter}. In this case, the order quantity can be viewed as a
weighted average of the ``textbook" order quantity $x_t^*$, and the
estimated mean $\hat{\mu}_t$. A higher $\gamma$ indicates a stronger PtC
effect.

Note that the model \eqref{eq:PtC}, like the DC model \eqref{eq:dc},
assumes that the adjustment parameter is constant over time and that
adjustments happen on each observation. We argue that these assumptions
are reasonable because they express a behaviour on average, similar to how
the DC behaviour is modelled via \eqref{eq:dc}.


\section{An Alternative Perspective and Model} \label{se:alt}

In this section, we argue that there are some positive aspects to
na\"{\i}ve adjustment. We also present a new adjustment model, which
allows one to perform demand-chasing and pull-to-centre in combination
if desired. 

\subsection{In favour of na\"{\i}ve adjustment}

As one can see from Section \ref{se:lit}, the previous literature on
judgemental adjustment has assumed, either implicitly or explicitly, that
DC and PtC are harmful. In this paper, we take a different point of view:
we argue that DC and PtC may sometimes be \emph{beneficial} in practice.
To see why, note that the ``textbook" NVP formula (\ref{eq:textbook})
applies only when one has an accurate statistical model of the demand
distribution. In reality, of course, such a model is rarely available. As
a result, the textbook formula may give the wrong answer in practice,
either underestimating or overestimating the optimal order level. In
some circumstances, therefore, DC and/or PtC might help rather than
hinder.

To be more specific, we suggest that a modest amount of ``na\"{\i}ve"
adjustment may be beneficial in two practically important situations:
\begin{enumerate}
\item When the demand model is correct, but there is insufficient data to
estimate its parameters accurately.
\item When the demand model is mis-specified.
\end{enumerate}
We will test these hypotheses using simulation experiments in the next
section, modelling the two situations.

There is another key difference between our work and the existing
literature. As mentioned above, the latter relies almost exclusively on
data collected from behavioural experiments with human subjects. Here, by
contrast, we will use simulated data, since it allows us to conduct
extensive experiments very easily.

\subsection{An integrated adjustment model}

To proceed, we make some additional remarks about the DC model
(\ref{eq:dc}). In our view, it is unlikely to be a good model of human
behaviour. Indeed, we have already observed that it effectively estimates the mean demand, and does not take cost information into account.

In an attempt to remedy the above weakness, we now propose a ``two-stage"
model of na\"{\i}ve adjustment, in which DC takes place after PtC:
\begin{equation} \label{eq:2stage}
    \begin{aligned}
        & x_t^\prime = (1-\gamma) x_t^* + \gamma \hat{\mu}_t \\
        & x_t = x_{t}^\prime + \beta \big( d_{t-1}-x_{t-1} \big).
    \end{aligned}
\end{equation}
The idea here is that we first take the ``textbook" order quantity
$x_t^*$, and apply PtC with parameter $\gamma$. This yields an adjusted
order quantity, here denoted by $x_t^\prime$. After that, we adjust
$x_t^\prime$ itself, by applying DC with parameter $\beta$.

Unlike the classical DC model (\ref{eq:dc}), the two-stage model
(\ref{eq:2stage}) yields non-trivial estimates of the optimal order
quantity, rather than merely estimating the mean demand. We remark that
we are not claiming that human practitioners actually use such a model
consciously. Indeed, it is likely that most decision makers will not make
a conscious distinction between the two effects. On the other hand, it
seems likely that, in practice, judgemental adjustments may contain
elements of both PtC and DC. In any case, our goal is not to explore the
behaviour of decision makers in practice, but to investigate whether the
model \eqref{eq:2stage} might be useful in the two specific situations
mentioned in the previous subsection.

Inserting the first equation in the second one in \eqref{eq:2stage}, we
obtain a unified formulation for DC and PtC, which summarises the order
adjustment in one formula:
\begin{equation} \label{eq:unified}
x_t = (1-\gamma) x^*_t + \gamma \hat{\mu}_t +
\beta \big( d_{t-1}-x_{t-1} \big).
\end{equation}
This makes it clearer that the adjusted order quantity is a linear
combination of three terms: system order quantity, mean and actual demand.
We will focus our investigation on the case in which the parameters
$\beta$ and $\gamma$ take values between 0 and 0.5 (This is a common
assumption in the literature, e.g. \cite{Be08,BK08,BHS08}), although the
theoretical parameter ranges might be wider.

It is important to note that, for a given $t$, the estimate $\hat{\mu}_t$
is itself based on $d_1, \ldots, d_{t-1}$, and so are the quantities
$x^*_t$ and $x_{t-1}$. Thus, all three quantities are subject to estimation
errors.

\section{Experiments on Simulated Data}
\label{se:exp}

In this section, we perform extensive computational experiments, to test the
two hypotheses as mentioned in the previous section. In Subsection
\ref{sub:method}, we describe our methodology. The hypotheses themselves are
tested in Subsections \ref{sub:known} and \ref{sub:unknown}, respectively.

We assume initially that $\tau = 0.7$, a value commonly used in the NVP
literature to approximate real-life problems (e.g., fashion retail, nurse
staffing) (\cite{AE05,LP01}). We show later in this section that the results
of our experiment hold with other NVP parameters as well.

\subsection{Methodology}
\label{sub:method}

The first step is to construct 500 time series, each consisting of 200
consecutive demand realisations, which is sufficiently long as shown in
behavioural experiments (\cite{Be08,SC00,BHS08}). To do this, we use the
{\tt arima.sim()} function from the {\tt stats} package in {\tt R}. We
assume that the ``true" DGP for the demands is an ARIMA(1,0,1) process, with
an initial mean of 10,000. We also assume that the noise term is normally
distributed with a standard deviation of 100.

We use an ARIMA model because it is popular in the NVP literature, and we
choose a model with two parameters so that we can explore the effects of
both over- and under-parametrisation. For a given time series and a given
$t = 1, \ldots, 200$, we let $d_t$ denote the demand realisation in the
$t^{th}$ time period.

Now suppose that we have selected a forecasting model. This can be the
correct model, i.e., ARIMA(1,0,1), or an incorrectly specified model,
such as AR(1). Suppose also that we have selected the adjustment parameters
$\beta$ and $\gamma$. We do the following for each time series:
\begin{enumerate}
\item For $t = 21, \ldots, 200$, we use the {\tt arima()} function in the
{\tt stats} package to produce maximum-likelihood estimates of the mean and
standard deviation of demand in time period $t$. We let $\hat{\mu}_t$ and
$\hat{\sigma}_t$ denote these estimates. We note that for $t \le 20$, the
results may be biased due to the shortage of observations.
\item For $t = 21, \ldots, 200$, we use $\tau$, $\hat{\mu}_t$ and
$\hat{\sigma}_t$ to compute the ``textbook" optimal order quantity for time
period $t$ using the formula (\ref{eq:textbook}). We let $x^*_t$ denote this
quantity.
\item Finally, we simulate the adjustment process with $\beta$ and $\gamma$.
To avoid systematic bias, we assume that $x_{20} = d_{20}$. For
$t = 21, \ldots, 200$, we assume that the amount ordered at the start of
period $t$ follows formula (\ref{eq:unified}).
\end{enumerate}

To quantify the effect of adjustment, we proceed as follows. For a given
series and for $t = 21, \ldots, 200$, we compute
\begin{equation}
\mbox{PPL} \big( x_t \big)  = 
100 \left[ \frac{\pi(d_t,d_t)-\pi(x_t,d_t)}{\pi(d_t,d_t)}
\right].
\end{equation}
(Here, ``PPL" stands for `percentage profit loss"; see \cite{LLS22}.) We
also compute the ``relative profit improvement" (also known as ``forecast
value added'' in some contexts (\cite{Gi10})):
\begin{equation}
   \mbox{RPI}(x_t) =1-\frac{PPL(x_t)}{PPL(x^*_t)}. 
\end{equation}
Intuitively, the mean of the RPI$(x_t)$, over all 500 time series,
represents the improvement in profit (if any) gained by using the chosen
adjustment in period $t$. The higher the value is, the better the
performance of the approach is. If the value is negative then this means
that the approach is worse than the benchmark.

\subsection{When the DGP is known}
\label{sub:known}

We first report results for the case in which the DGP is known, but the model
parameters need to be estimated. In particular, we assume that we are using
an ARIMA(1,0,1) model, but with unknown parameters.

In Table \ref{tab:compare_1/2_high}, we show the mean RPI for different values
of the adjustment parameters. The heading ``short dataset" indicates that the
mean RPI is computed over the interval $t \in [21,110]$, and the heading
``long dataset" indicates that the mean is computed over the interval
$t \in [111,200]$. 

\begin{table}[htb]
\caption{Average RPI with varying adjustment parameters and with
short/long datasets ($\tau =0.7$)} \label{tab:compare_1/2_high}
\centering
\resizebox{0.7\linewidth}{!}{
\begin{tabular}{ccccccc}
\toprule
\multicolumn{1}{c}{\textbf{Short dataset}} & \multicolumn{6}{c}{$\gamma$}\\
\cmidrule(l{3pt}r{3pt}){1-1} \cmidrule(l{3pt}r{3pt}){2-7}
$\beta$ & $0$ & $0.1$ & $0.2$ & $0.3$ & $0.4$ & $0.5$ \\
\midrule
$0$ & $-$ & $3.5\%$ & $3.4\%$ & $3.3\%$ & $3.2\%$ & $2.8\%$\\
$0.1$ & $3.8\%$ & $\textbf{5.0\%}$ & $4.8\%$ & $4.7\%$ & $4.6\%$ & $3.8\%$\\
$0.2$ & $4.2\%$ & $\textbf{5.1\%}$ & $\textbf{5.0\%}$ & $4.6\%$ & $3.0\%$ & $2.4\%$\\
$0.3$ & $2.5\%$ & $2.8\%$ & $2.6\%$ & $1.5\%$ & $0.9\%$ & $0.2\%$\\
$0.4$ & $-2.4\%$ & $-1.3\%$ & $-1.9\%$ & $-3.1\%$ & $-3.9\%$ & $-5.0\%$\\
$0.5$ & $-5.1\%$ & $-3.9\%$ & $-4.4\%$ & $-5.9\%$ & $-6.8\%$ & $-7.7\%$\\
\midrule
\multicolumn{1}{c}{\textbf{Long dataset}} & \multicolumn{6}{c}{$\gamma$}\\
\cmidrule(l{3pt}r{3pt}){1-1} \cmidrule(l{3pt}r{3pt}){2-7}
$\beta$ & $0$ & $0.1$ & $0.2$ & $0.3$ & $0.4$ & $0.5$ \\
\midrule
$0$ & $-$ & $-0.5\%$ & $-0.9\%$ & $-1.3\%$ & $-1.7\%$ & $-2.1\%$\\
$0.1$ & $-0.7\%$ & $-0.2\%$ & $-0.5\%$ & $-0.7\%$ & $-0.9\%$ & $-1.3\%$\\
$0.2$ & $0.1\%$ & $0.3\%$ & $0.2\%$ & $0.1\%$ & $-0.2\%$ & $-0.7\%$\\
$0.3$ & $-2.1\%$ & $-1.6\%$ & $-2.6\%$ & $-3.3\%$ & $-4.1\%$ & $-5.1\%$\\
$0.4$ & $-6.7\%$ & $-6.4\%$ & $-7.5\%$ & $-8.4\%$ & $-9.2\%$ & $-10.2\%$\\
$0.5$ & $-11.1\%$ & $-10.5\%$ & $-11.8\%$ & $-12.8\%$ & $-13.8\%$ & $-14.7\%$\\
\bottomrule
\end{tabular}}
\end{table}

Table \ref{tab:compare_1/2_high} indicates that a modest amount of adjustment
can be beneficial, especially when the number of historical demand
observations is short. This is possibly due to the statistical model being
unable to estimate its parameters accurately on insufficient data. Therefore,
modest amount of adjustment can provide additional information. On the other
hand, too much adjustment leads to a loss. We also find that the RPI is more
sensitive to the choices of $\beta$, meaning that demand chasing has a bigger
influence in this case. We mark that this is also true for other data inputs.

\begin{figure}[!htb]
 \centering
 \caption{Average RPI of judgemental adjustments vs. dataset length
 ($\tau = 0.7$)} \label{fig:RPI vs length}
 \includegraphics[width=0.85\textwidth]{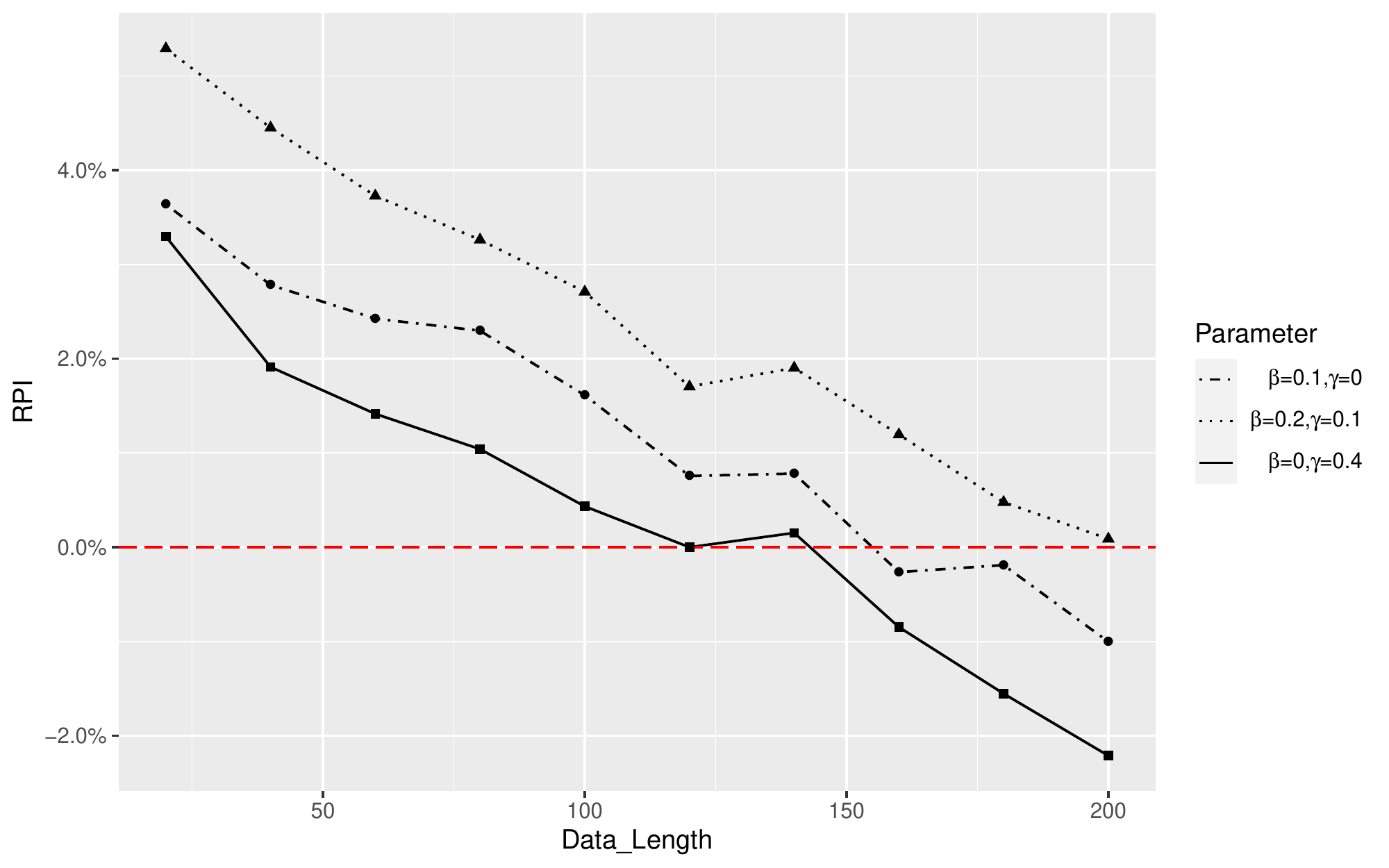}
\end{figure}

To explore this effect in more detail, we show in Figure
\ref{fig:RPI vs length} a plot of average RPI against the length of the
dataset, for three different values of $(\beta,\gamma)$, namely, $(0,0.4)$,
$(0.1,0)$ and $(0.2,0.1)$.  It can be seen that the average RPI is well
above zero initially, but decreases, and eventually becomes negative.

A tentative explanation is that the maximum-likelihood estimates of
$\hat{\mu}_t$ and $\hat{\sigma}_t$ are prone to errors when the number of
observations is small. It may even be that the estimates suffer from some
kind of systematic bias, which decreases over time. By performing a small
amount of adjustment, we shift the order quantity toward the true
optimal value. We would expect this effect to vanish as more
data becomes available.

\begin{figure}[htb]
 \centering
 \caption{Heat map of RPIs for different combinations of adjustment parameters
 (data length = 20)}
 \label{fig:RPI heat}
 \includegraphics[width=0.85\textwidth]{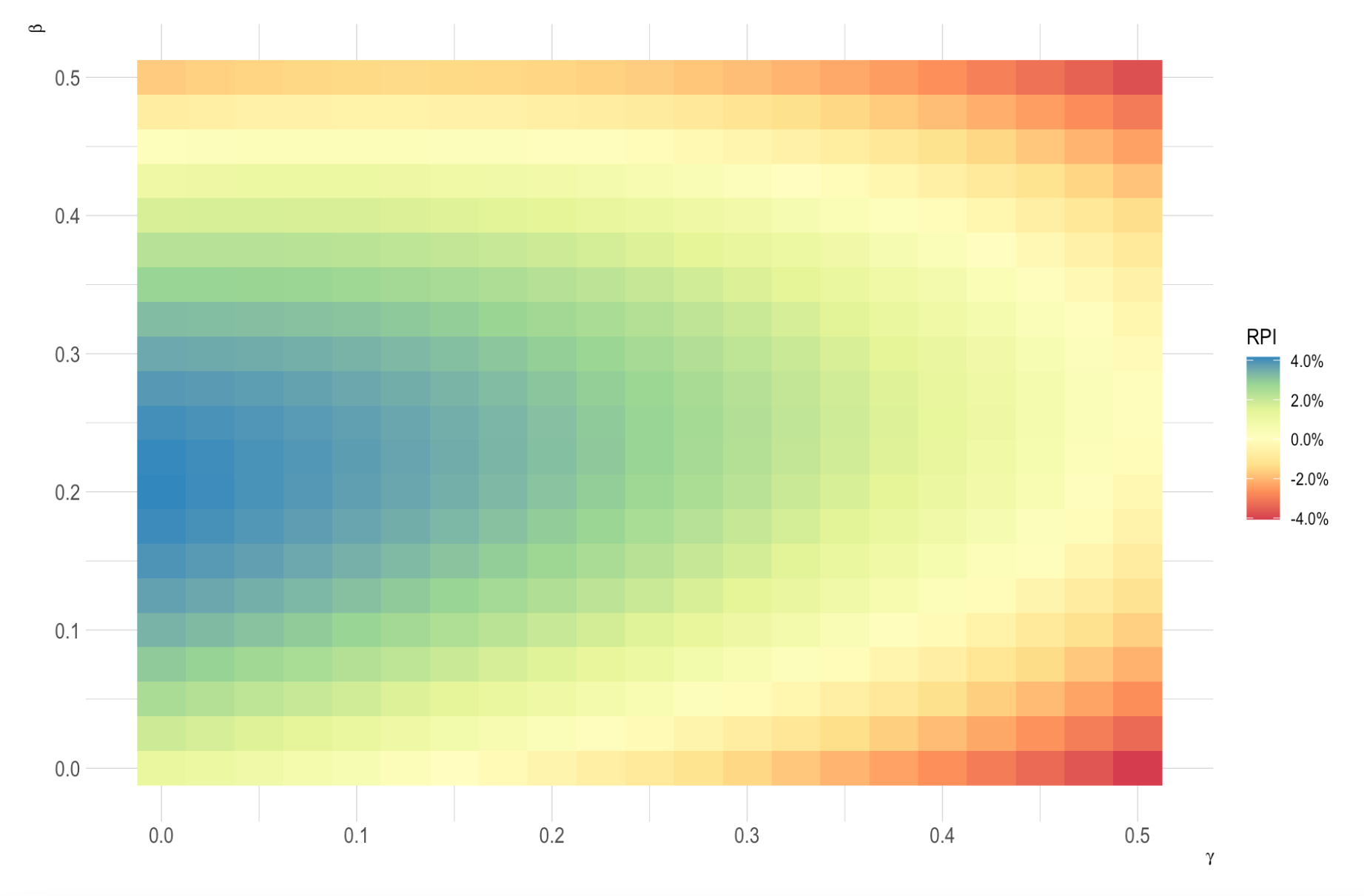}
\end{figure}

In Figure \ref{fig:RPI heat}, we use a heat map to compare the performance
of adjustment with different parameter values, for the case in which the
number of past observations is exactly 20 (i.e., $t=21$). As we observed in
Figure \ref{fig:RPI vs length}, the dominance relationship between parameter
pairs (one parameter pair generating higher RPI than the other) is not
influenced by data length. Therefore, the choice of $t=21$ can amplify the
results, making it easier for us to observe. One can clearly see from
Figure \ref{fig:RPI heat} that large adjustments are harmful, while modest
amounts of adjustment, on the other hand, can generate positive RPI. From
our data, the most profitable option is to set $\beta$ to around $0.2$ and
$\gamma$ to around $0.05$. A possible explanation to this phenomenon is that
the ``textbook'' order quantity will be close to the theoretical one, but
it would need to be adjusted by the order quantity and demand on the
previous observation. The PtC effect needs to be reduced in comparison with
the DC effect.

\begin{figure}[htb]
 \centering
 \caption{Boxplots of the RPI for different values of $\tau$ and different
 combinations of adjustment parameters}
 \label{fig:RPI vs tau}
 \includegraphics[width=0.85\textwidth]{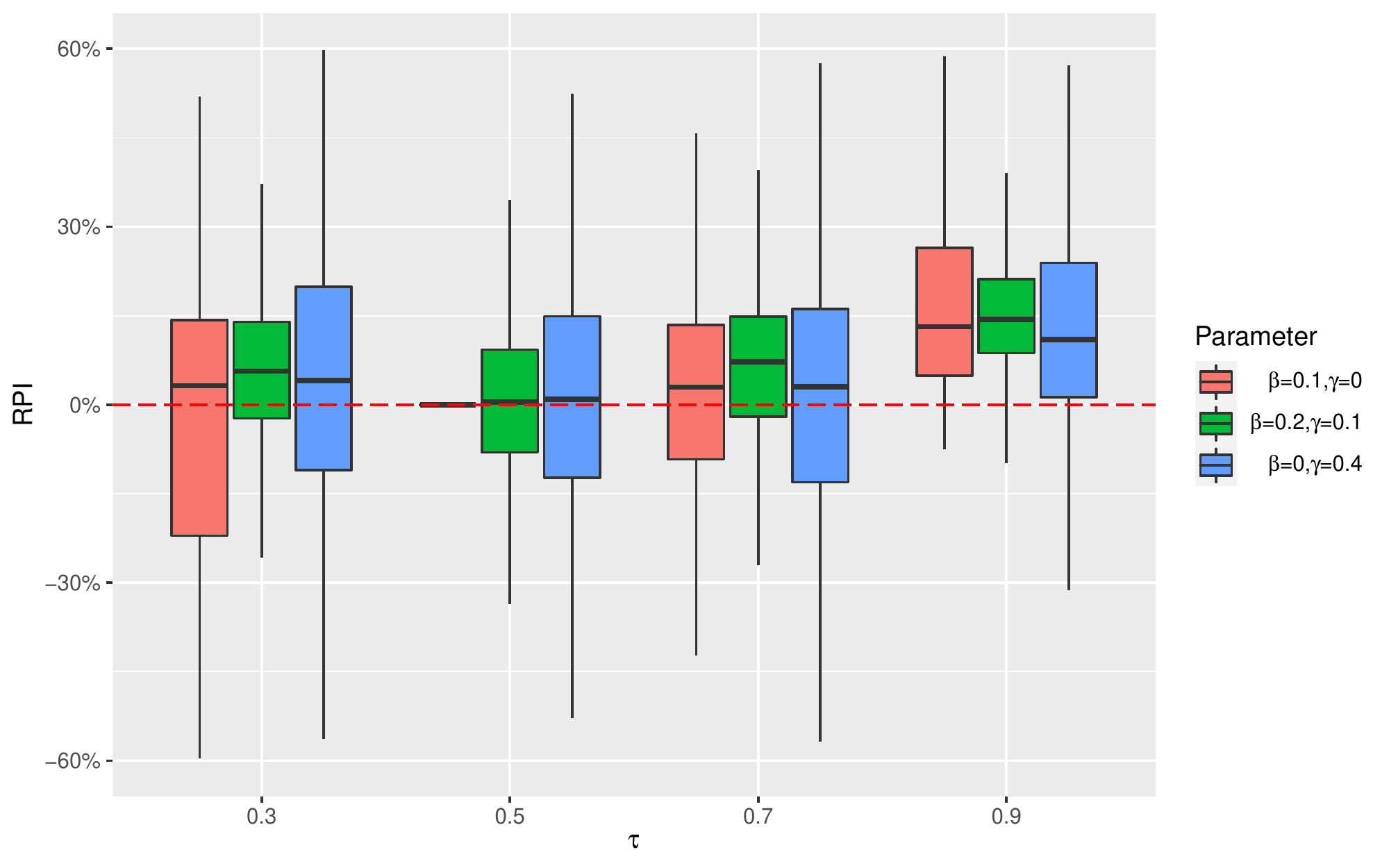}
\end{figure}

Next, we examine the effect of $\tau$ on the performance of adjustment.
Figure \ref{fig:RPI vs tau} shows boxplots of the RPI for four different
values of $\tau$ ($0.3$, $0.5$, $0.7$ and $0.9$), and the same three
values for $(\beta,\gamma)$ as before. Here, $t=21$ as before. Each boxplot
shows the range, median and quartiles over the 500 time series.

Interestingly, adjustment yields a benefit in every case, \emph{except}
when $\tau=0.5$. Moreover, we can see in Figure \ref{fig:RPI vs tau} that
the performance of adjustments when $\tau = 0.3$ is very similar to that
when $\tau = 0.7$. This suggests that the effect of adjustment may be
symmetric around $0.5$. Note also that, even when $\tau=0.5$, adjustment
does not cause any noticeable loss of profit.

All things considered, it appears that, when the model is correct but the
data length is short, a modest amount of na\"{\i}ve adjustment can be
beneficial instead of harmful. This goes against the prevailing view in
the literature that DC and PtC are invariably damaging. We believe that
the discrepancy is mainly due to the model correctness assumption made in
the literature. When model correctness is assured, there is no doubt that
any kind of adjustment will be harmful to the profit. Indeed, the
behavioural results are in line with the results that we obtained with the
longest datasets in our experiment, where there was sufficient data to
properly estimate all parameters.


\subsection{When the model is misspecified}
\label{sub:unknown}

In this subsection, we examine the effect of model misspecification on the
relative performance of the judgemental adjustments. We consider three
scenarios of model misspecification:
\begin{enumerate}
 \item The model is under-parametrised (i.e., omits one or more important
 variables), which typically leads to biased estimates of parameters;
 \item The model is over-parametrised (i.e., has one or more redundant
 variables), which usually leads to inefficient estimates of parameters;
 \item The model has the correct parameters, but the assumed distribution of
 the error term is wrong, which can lead to biased quantile estimates.
\end{enumerate}

For scenario \#1, we use an MA(1) model to fit the underlying demand data.
For scenario \#2, we use an ARIMA(2,0,1) model. For scenario \#3, the data
is generated using a modified version of ARIMA model, in which the error
term follows the Laplace distribution instead of the normal distribution.
When estimating the optimal order quantity, however, we use the incorrect
assumption that the error term follows the normal distribution.

Since we wish to focus on the effect of model misspecification, rather than
the effect of a lack of data (as in the previous subsection), we report the
mean RPIs when $t=200$, when plenty of demand data is available. As before,
however, all means are taken over 500 time series.

\begin{figure}[ht]
 \centering
 \caption{Boxplots of the RPI for three different models and three different
 combinations of adjustment parameters ($\tau=0.7$, $t=200$)}
 \label{fig:RPI mis}
 \includegraphics[width=0.85\textwidth]{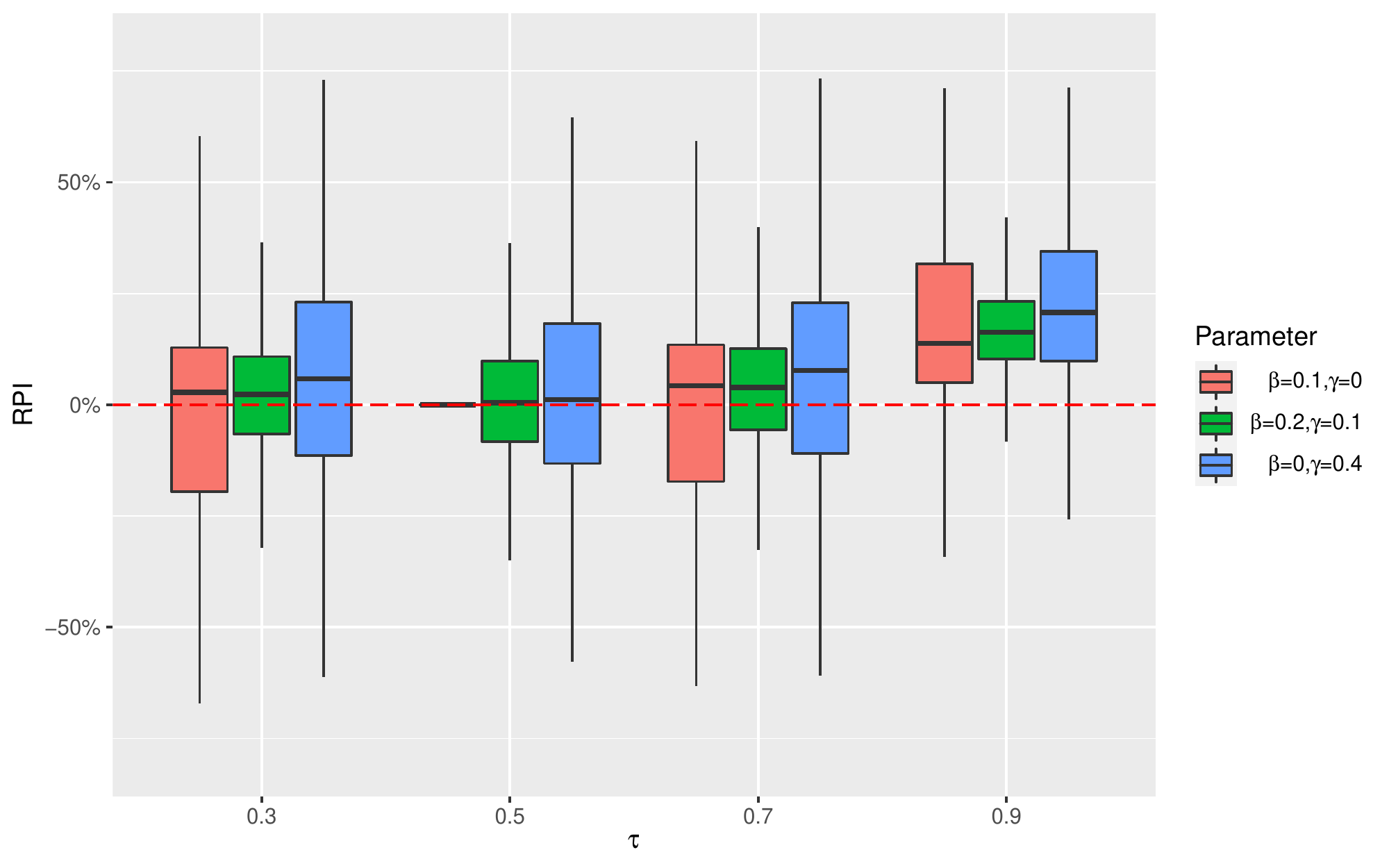}
\end{figure}

Figure \ref{fig:RPI mis} shows the boxplots for scenarios \#1 and scenario
\#2, together with the boxplots for the correctly specified case, for
comparison. Here, $\tau$ is equal to $0.7$. As before, results are reported
for three different settings of $(\beta,\gamma)$.

It is apparent that, when the forecasting model omits important variables,
na\"{\i}ve adjustment can yield significant increase in profit. The
improvement seems to be present also when the model has redundant variables,
but the effect is less pronounced. This is probably due to the nature of
judgemental adjustment. In the case where important information is omitted
from the model, there is a good chance that adjustment can provide
supplementary information. On the other hand, when the model
contains redundant variables, no additional information is needed.

Heatmaps for scenario \#1 and scenario \#2 are presented in Appendix
\ref{app:heat}. They confirm again our findings in Figure \ref{fig:RPI mis}.
It appears that the optimal PtC parameter ($\gamma$) is larger than the
optimal DC parameter ($\beta$) in the under-parametrised case. This is
probably because omitting variables in the model leads to an over-estimate
of the variance, leading to order quantities that are further from the mean
than required. This bias is remedied by using PtC, since the order quantity
is brought closer to the mean of the data. On the other hand, the
improvement from adjustment when the model has redundant variables is
limited, because the model overfits the data, thus having a variance biased towards zero.

\begin{figure}[ht]
 \centering
 \caption{Boxplots of the RPI when assumed distribution of error term
 is wrong ($t=200$)} \label{fig:RPI error}
 \includegraphics[width=0.85\textwidth]{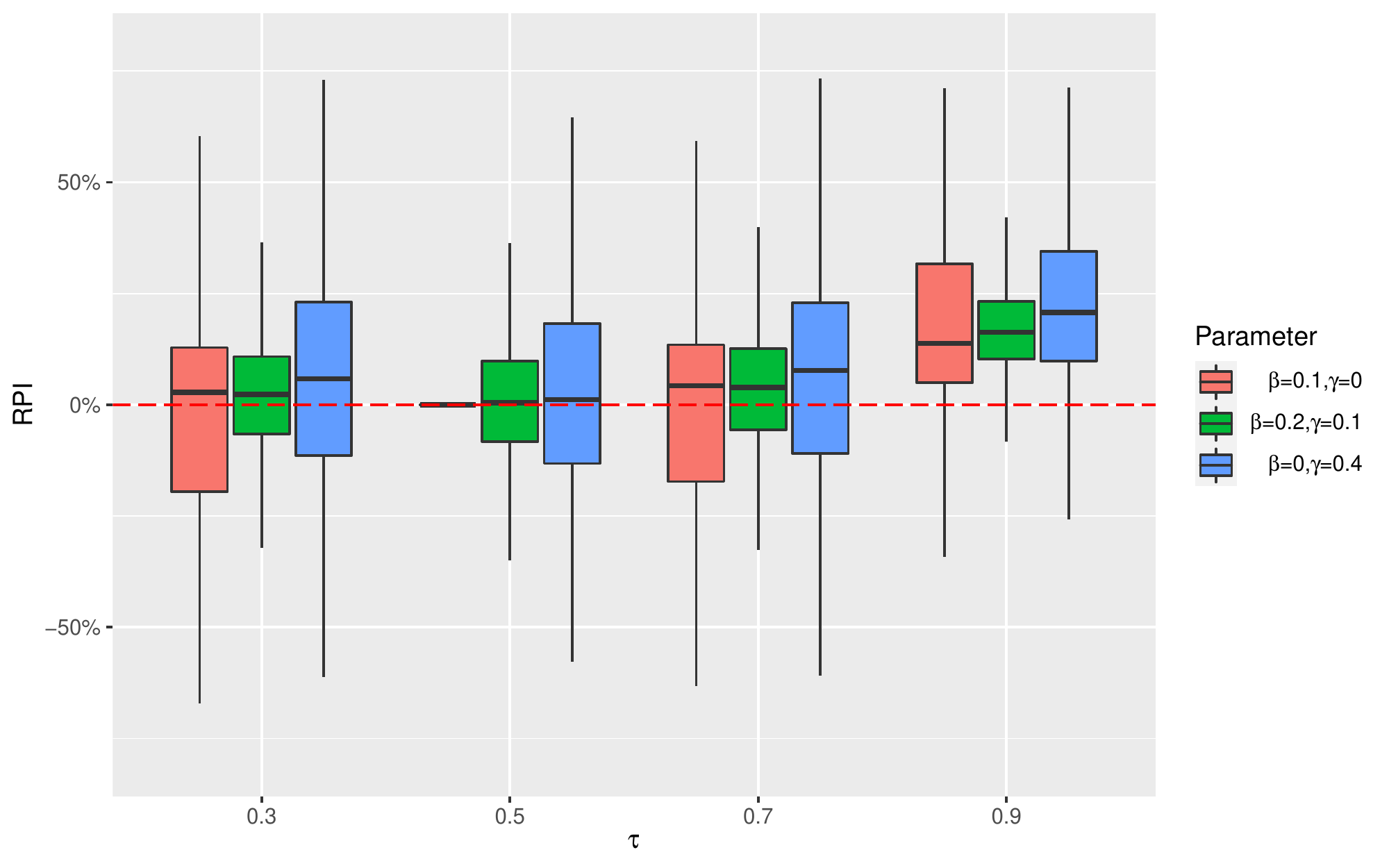}
\end{figure}

Finally, we consider scenario \#3, in which the underlying assumption on the
error term distribution is wrong. Figure \ref{fig:RPI error} shows the RPI
boxplots for this case. It can be seen that adjustment significantly
improves the expected profit when $\tau=0.9$, but the effect disappears
when $\tau = 0.5$. This can be explained by the relative shapes of the
normal and Laplace distributions. Since the largest difference between the
distributions is in the ``tails", more adjustment will be needed when
$\tau$ is at an extreme value (i.e., close to either $0$ or $1$).

Interestingly, PtC seems to be of more benefit than DC in scenario \#3.
This is probably  because the wrong assumption about the error term
distribution leads to a systematic over-estimation of the variance. Just as
in the under-parametrised case, this over-estimation alleviated by PtC.

Once again, we find that na\"{\i}ve adjustments can be beneficial instead
of harmful. The results are similar to what was found in Subsection
\ref{sub:known} and can be overall summarised as:
\begin{enumerate}
    \item DC is more useful in the short data length case. This is probably
    because of the small sample bias, which vanishes on larger samples.
    \item PtC is more beneficial than DC in the case of model
    misspecification. This is because of the systematic underestimation of
    variance in case of omitted variables or wrong model form.
    \item In the case when the DGP is known, the benefits from adjustments
    only depend on data length and adjustment parameters.
    \item The improvement from adjustments appears when the model omits
    important variables, or the distributional assumption is wrong.
    \item When the model has redundant variables, the improvement from
    both DC and PtC is less pronounced.
\end{enumerate}
In general, we can conclude that DC is likely to be beneficial when there
is a shortage of demand data, whereas PtC is likely to bring value when the
model is misspecified. The combination of the two should be beneficial in practice where the model is not known and the samples of data are limited.


\section{Tuning Algorithm}
\label{se:algorithm}

In the previous section, we showed that na\"{\i}ve adjustments can improve
the expected profit when the data is insufficient and/or the demand model
is misspecified. It is however not clear how one might choose suitable
parameter values when faced with a specific NVP instance. In this section,
we propose and test a simple heuristic algorithm for parameter ``tuning". We
believe that this tuning algorithm may be of interest to both academics and
practitioners. The method is explained in Subsection \ref{sub:pro}. In
Subsection \ref{sub:exam}, we apply our method to a real-life NVP instance,
for which the true model is not known. Finally, in Subsection
\ref{sub:result}, we present and compare the results with and without the
application of adjustment with ``tuned" parameters.

\subsection{Procedure} \label{sub:pro}

Let us suppose that we have a forecasting model that, for each period in the
training set, $t \in [1,s]$, is able to provide an estimate of the mean demand
$\hat{\mu}_t$. Moreover, let us assume that we are also able to estimate the
``textbook” optimal order quantity $x^*_t$. We remind the reader that our
adjusted order quantity takes the form:
\[
    x_t = (1-\gamma) x_t^* + \gamma \hat{\mu}_t + \beta (d_{t-1}-x_{t-1}).
\]
The $x_t$ can be calculated for all $t \in [1, s]$ based on the available mean and actual demand, the ``textbook'' order quantity and some values of $\beta$ and $\gamma$. To determine the values of parameters, we solve an optimisation problem over
the training set. Following suggestions in \cite{BR19} and \cite{LLS22}, we
use a non-standard ``loss function" for this purpose. The function is chosen
to maximise the profit over the training set, instead of minimising the MSE
or MAE in the usual way. That is, we estimate $\beta$ and $\gamma$ by
maximising the in-sample empirical profit:
\begin{equation}
\label{eq:max}
    \max_{\beta,\gamma} \sum_{t=2}^s {\pi(x_t^*,\hat{\mu}_t,x_t,d_t)}.
\end{equation}

We remark that the function to be maximised in (\ref{eq:max}) is continuous
and concave. On the other hand, it is not differentiable in general. This means that in order to maximise the profit we need to use derivative-free optimisation algorithms, such as Nelder-Mead (\cite{NM65,RS13}).

\subsection{Real life example} \label{sub:exam}

Here we present an example of application of our approach to real data.
The data we use comes from a medium-sized grocery store which sells a wide
range of products, many of which are perishable. It includes daily demands
for each product, for a period of around 9 weeks, which ran from mid-October
to December. For this study, we selected four typical products with very
different data structures and NVP parameters. In particular, we made sure
that the selected products have a range of prices, demands and critical
quantiles, in order to make the experiment less biased. For reasons of
confidentiality, we refer to these products as simply A, B, C and D.

\begin{figure}[htb]
    \centering
    \subfloat[Demand for product A]{%
    \resizebox*{5.2cm}{!}{\includegraphics{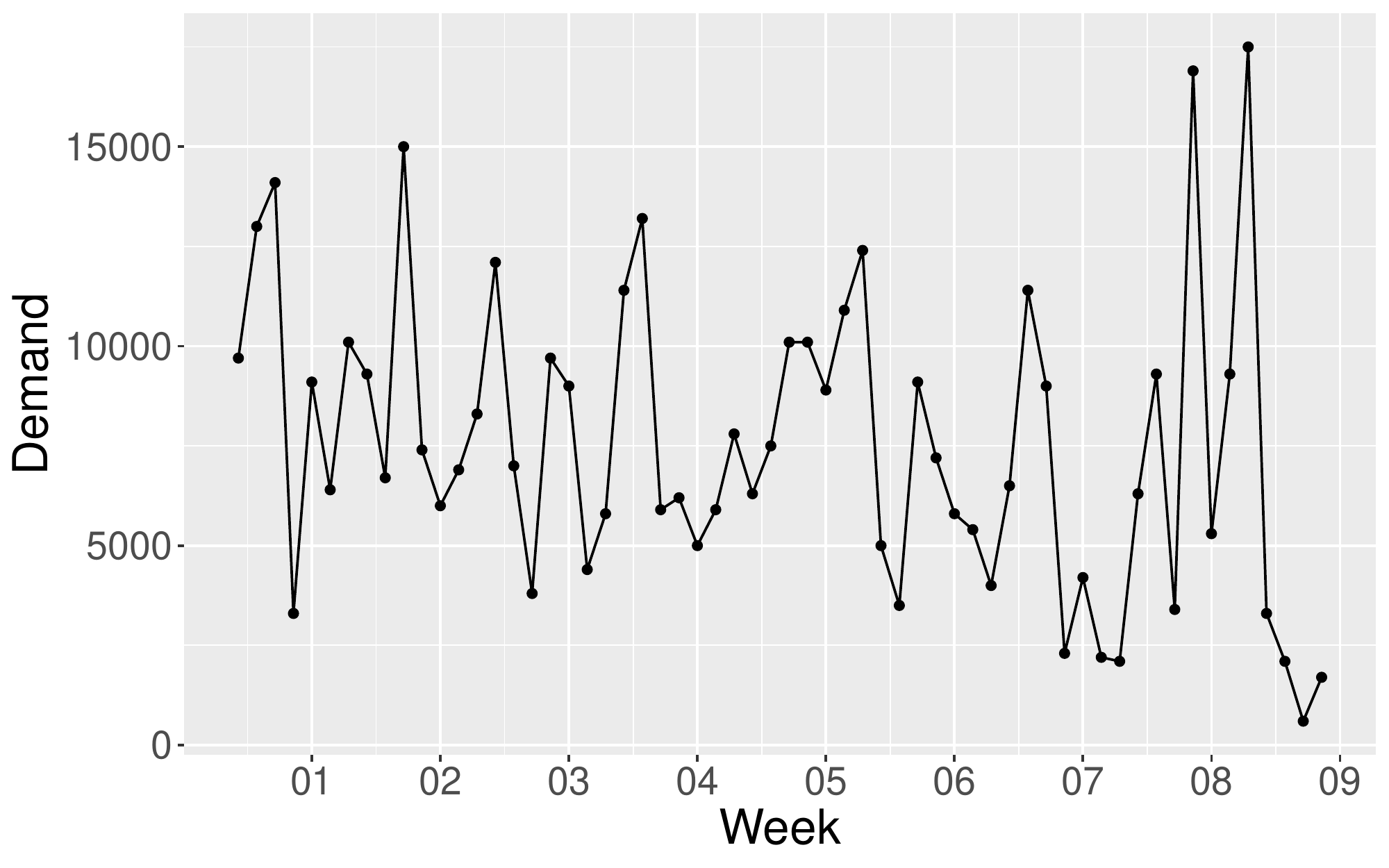}}}\hspace{5pt}
    \subfloat[Demand for product B]{%
    \resizebox*{5.2cm}{!}{\includegraphics{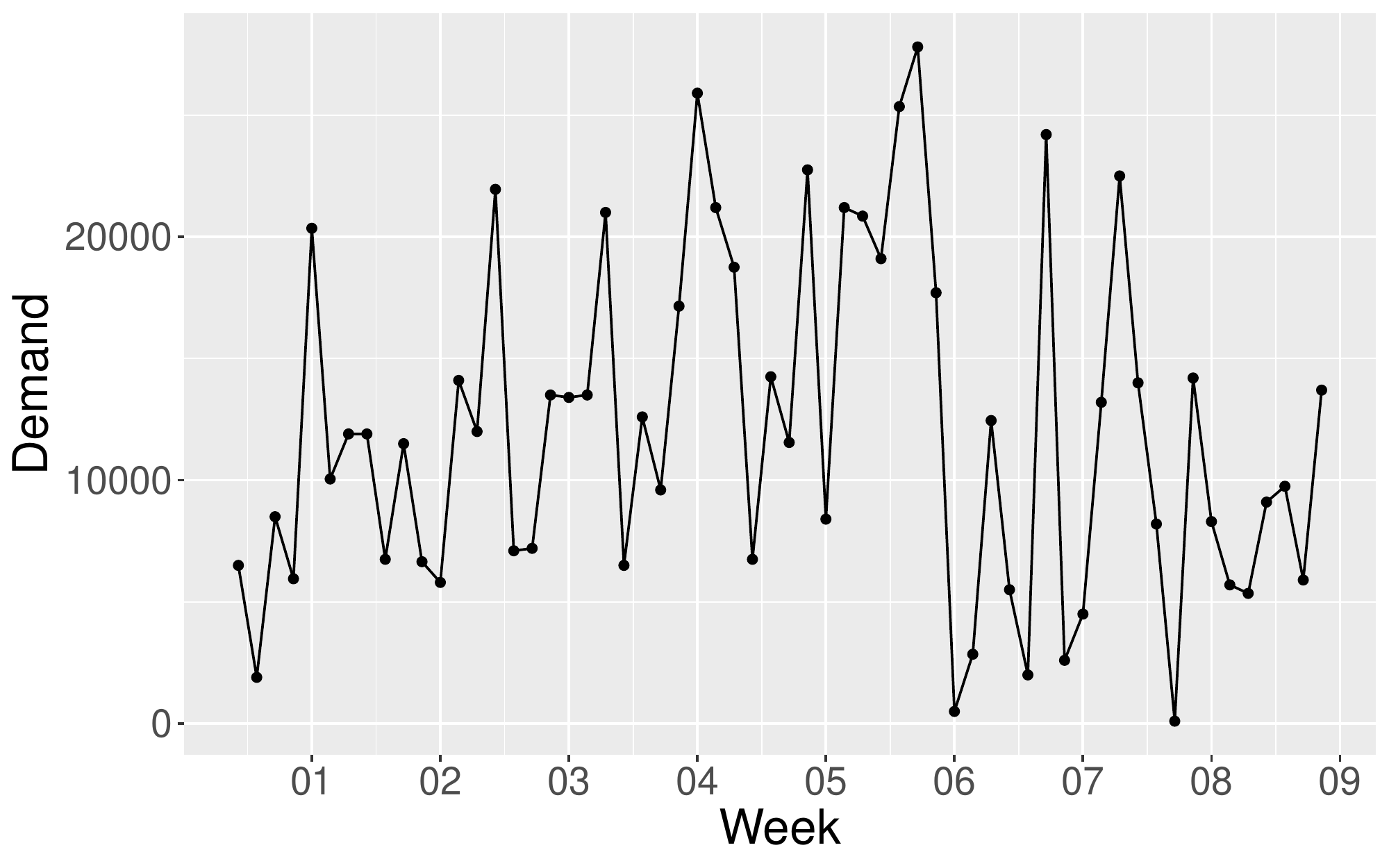}}}\\
    \subfloat[Demand for product C]{%
    \resizebox*{5.2cm}{!}{\includegraphics{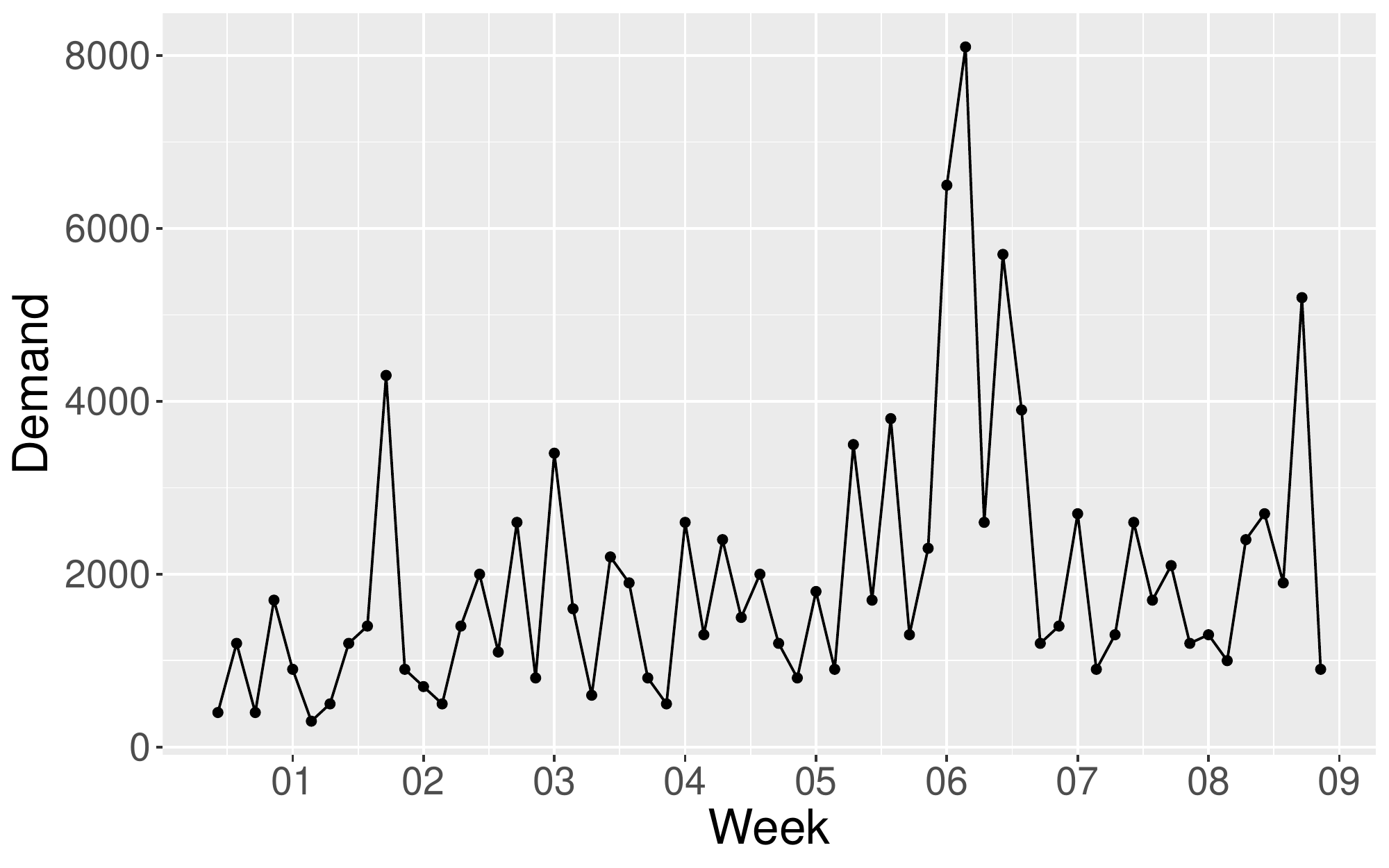}}}\hspace{5pt}
    \subfloat[Demand for product D]{%
    \resizebox*{5.2cm}{!}{\includegraphics{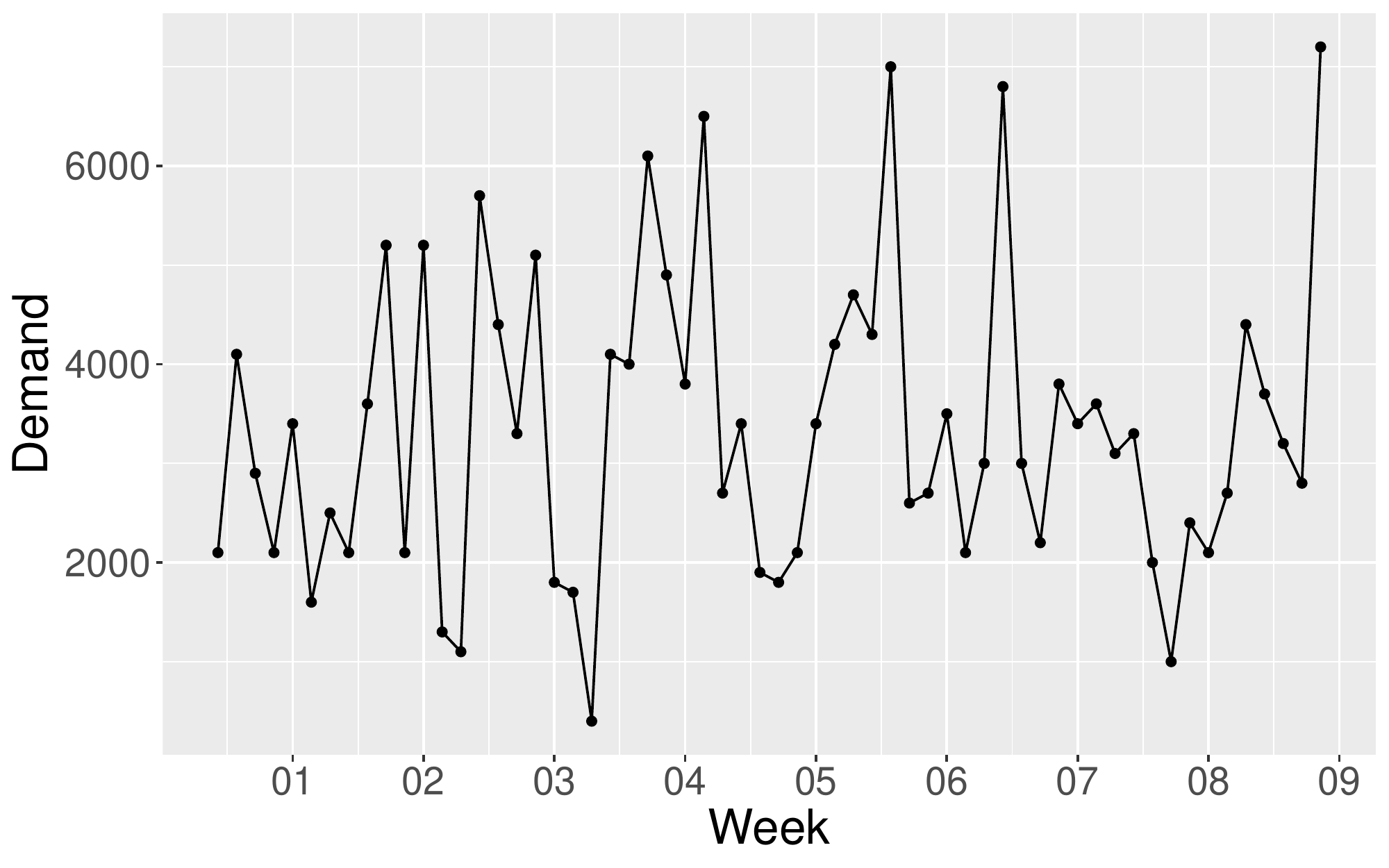}}}
    \caption{Demand time-series for real-life case}
    \label{fig:demand}
\end{figure}

To give reader some sense of the data, we provide time-series plots
in Figure \ref{fig:demand} and summarise the cost parameters in Table
\ref{tab:NVP_para}.

\begin{table}[htb]
\caption{Data for a subset of the products}
\label{tab:NVP_para}
\centering
\resizebox{0.64\linewidth}{!}{
\begin{tabular}{cccccc}
\toprule
\multicolumn{1}{c}{\textbf{Products}} & \multicolumn{4}{c}{\textbf{Price and Costs}} & \multicolumn{1}{c}{\textbf{Critical Quantile}} \\
\cmidrule(l{3pt}r{3pt}){1-1} \cmidrule(l{3pt}r{3pt}){2-5}
\cmidrule(l{3pt}r{3pt}){6-6}
& $p$ & $v$ & $c_h$ & $c_s$ &\\
\midrule
$A$ & 2.96 & 1.28 & 0.49 & 0.51 & 0.55\\
$B$ & 11.98 & 4.13 & 2.49 & 1.33 & 0.58\\
$C$ & 2.86 & 1.96 & 0.78 & 0.56 & 0.35\\
$D$ & 4.29 & 3.24 & 1.03 & 0.21 & 0.23\\
\bottomrule
\end{tabular}}
\end{table}

Following standard practice in forecasting, we use a rolling-origin method (\cite{Ta00}),
with constant in-sample size. For each product, on each iteration, we use
three-fifths of the data as the training set, and perform a one-step-ahead
forecast.

To perform a fair comparison, and reduce the possibility of bias in our
choice of model, we simply applied one of the most popular automatic
techniques for forecasting: the \texttt{ets()} function from the R
{\tt forecast} package (\cite{Hy20}). This function attempts to select the
most appropriate ETS model, using the Akaike Information Criterion. The
pre-tuning decision ($x_t^*$) is computed using the output from the
traditional forecasting procedure, while the tuned decision ($x_t$) is
computed using the output from the forecasting and tuning procedures in
combination.

\subsection{Results} \label{sub:result}

\begin{figure}[htb]
    \centering
    \caption{Boxplot of the out-of-sample RPI. The black lines in the boxes represent mean values}
    \label{fig:PPL real}
    \includegraphics[width=0.8\textwidth]{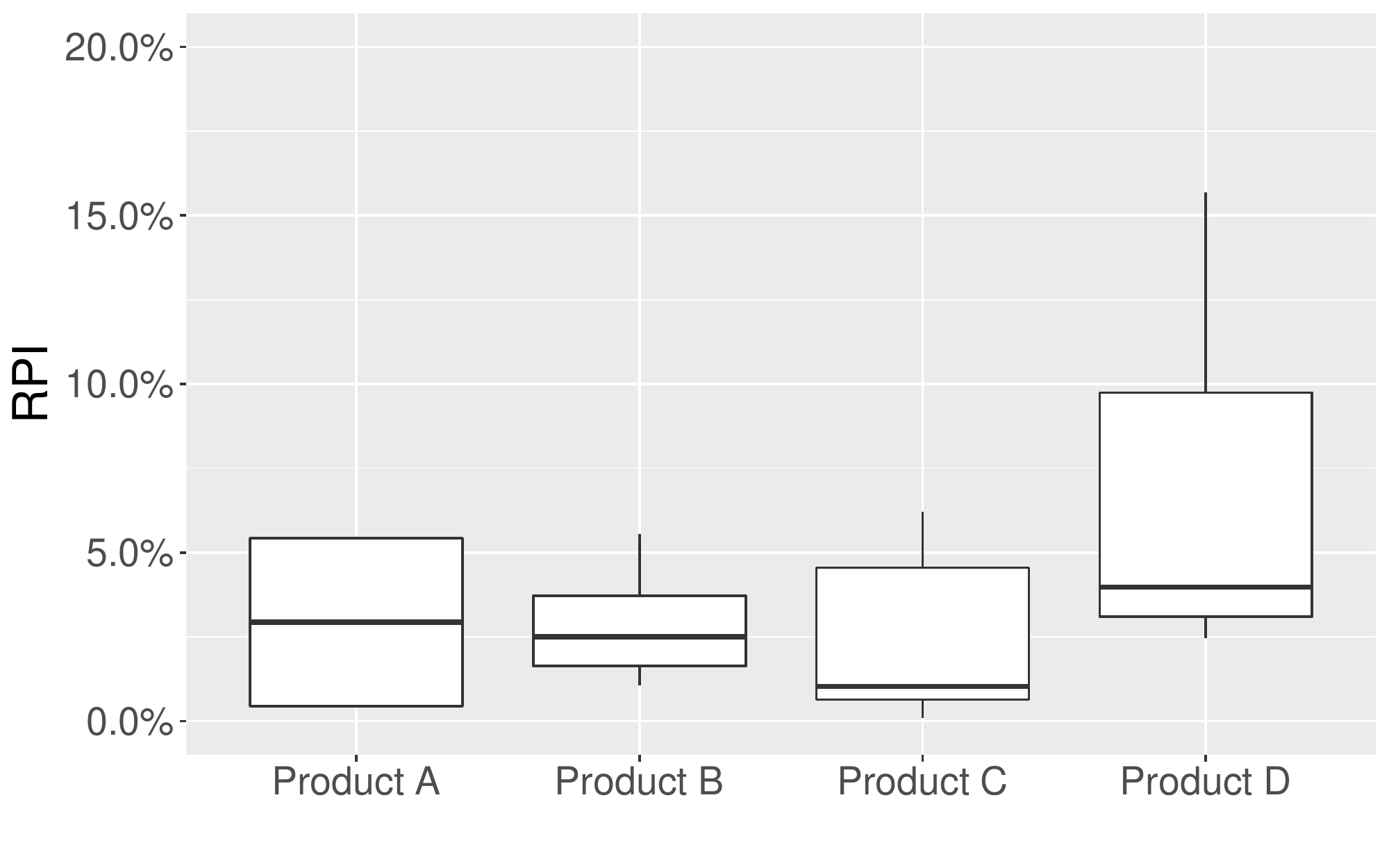}
\end{figure}

We now present the results obtained with our tuning algorithm. Figure
\ref{fig:PPL real} displays box-plots of the RPI, taken over the iterations,
for each of the four products. We remind the reader that a positive RPI
indicates that adjustment has been beneficial (see Subsection
\ref{sub:method}). The plots indicate that the RPI is positive for all
four products in all the situations. Thus, the tuned order
decisions outperform the pre-tuned ones for all four products.

Table \ref{tab:service} summarises performance in terms of service levels.
The row labelled `target' shows the critical quantile that maximises the
expected profit for the given cost and price parameters. The next two
rows show the achieved service level without and with tuning,
respectively. In all four cases, the achieved service level with tuning
is much closer to the target one than the one without tuning.

\begin{table}[htb]
\centering
\caption{Achieved Service level of each methods}
\label{tab:service}
\begin{tabular}{ccccc}
\toprule
& Product A & Product B & Product C & Product D\\
\midrule
Target & 0.55 & 0.58 & 0.35 & 0.23\\
Pre-tuning & 0.75 & 0.80 & 0.13 & 0.17\\
Tuned & 0.65 & 0.62 & 0.25 & 0.29\\
\bottomrule
\end{tabular}
\end{table}

To gain additional insight, we repeated the entire experiment using four
other popular forecasting methods. The results are presented in Table
\ref{tab:extensive}. From the table, one can see that the tuning algorithm
yields a positive out-of-sample RPI in every single case.

\begin{table}[htb]
\caption{RPI results obtained when using adjustment with other forecasting
methods}
\label{tab:extensive}
\centering
\resizebox{0.6\linewidth}{!}{
\begin{tabular}{c|ccc}
\toprule
\textbf{Products} & \textbf{Methods} & \textbf{In-sample} & \textbf{Out-of-sample}\\
\midrule
Product A & Mean & 4.2\% & 2.8\%\\
& S-Mean & 0.4\% & 2.4\%\\
& S-Naïve & 3.0\% & 2.4\%\\
& ARIMA & 1.0\% & 0.4\%\\
Product B & Mean & 3.9\% & 4.1\%\\
& S-Mean & 1.1\% & 2.9\%\\
& S-Naïve & 2.5\% & 2.8\%\\
& ARIMA & 1.9\% & 1.2\%\\
Product C & Mean & 4.1\% & 4.2\%\\
& S-Mean & 2.8\% & 2.2\%\\
& S-Naïve & 2.6\% & 2.5\%\\
& ARIMA & 0.2\% & 1.1\%\\
Product D & Mean & 4.8\% & 3.1\%\\
& S-Mean & 1.2\% & 3.5\%\\
& S-Naïve & 3.3\% & 4.5\%\\
& ARIMA & 2.6\% & 1.1\%\\
\bottomrule
\end{tabular}}
\end{table}

The experiment in this section supports our findings in the simulation study. We show that because the true model is not known, the proposed tuning algorithm leads to improvements, bringing the order closer to the correct level.


\section{Concluding Remarks} \label{se:end}

Although there is a considerable literature on judgemental adjustment for
newsvendor problems, it has been assumed up to now that `demand chasing'
and `pull-to-centre' are especially na\"{\i}ve, and likely to lead to
losses in profit. In this paper we have shown that, surprisingly, these
adjustment procedures can lead to increased profits in some situations.
In particular, they can be useful when (a) there is not enough data
available to estimate parameters accurately, and (b) the demand model is
misspecified. Interestingly, DC appears to be more useful under condition
(a), while PtC seems to be of more benefit under condition (b). In general,
this is because in case of (a) the order estimates suffer from some kind of
systematic bias due to short data length; while in the situation (b) the
estimated variance is often higher than needed.

We also proposed a simple heuristic for tuning the adjustment parameters.
Using a real-life example, we show that the tuned orders outperform the
pre-tuned ones in terms of the achieved profit, and also led to a service
level closer to the target one.

There are several interesting topics for further research. First, one could
attempt to characterise other scenarios under which na\"{\i}ve adjustments
tend to be beneficial. Second, one could examine the effects of other forms
of adjustment. Third, it might be beneficial to conduct behavioural
experiments, in the lab and/or field, to confirm the simulation results.
Finally, it would be interesting to extend the research to multi-item
newsvendor problems, either with or without substitution effects between
products.


\bibliographystyle{plain}
\bibliography{newsvendor3}


\newpage
\appendix

\section{Heatmaps for model misspecification} \label{app:heat}
\begin{figure}[h]
 \centering
 \caption{RPI heat map for the under-parametrised case
 ($\tau = 0.7$, $t=200$)}
 \includegraphics[width=0.8\textwidth]{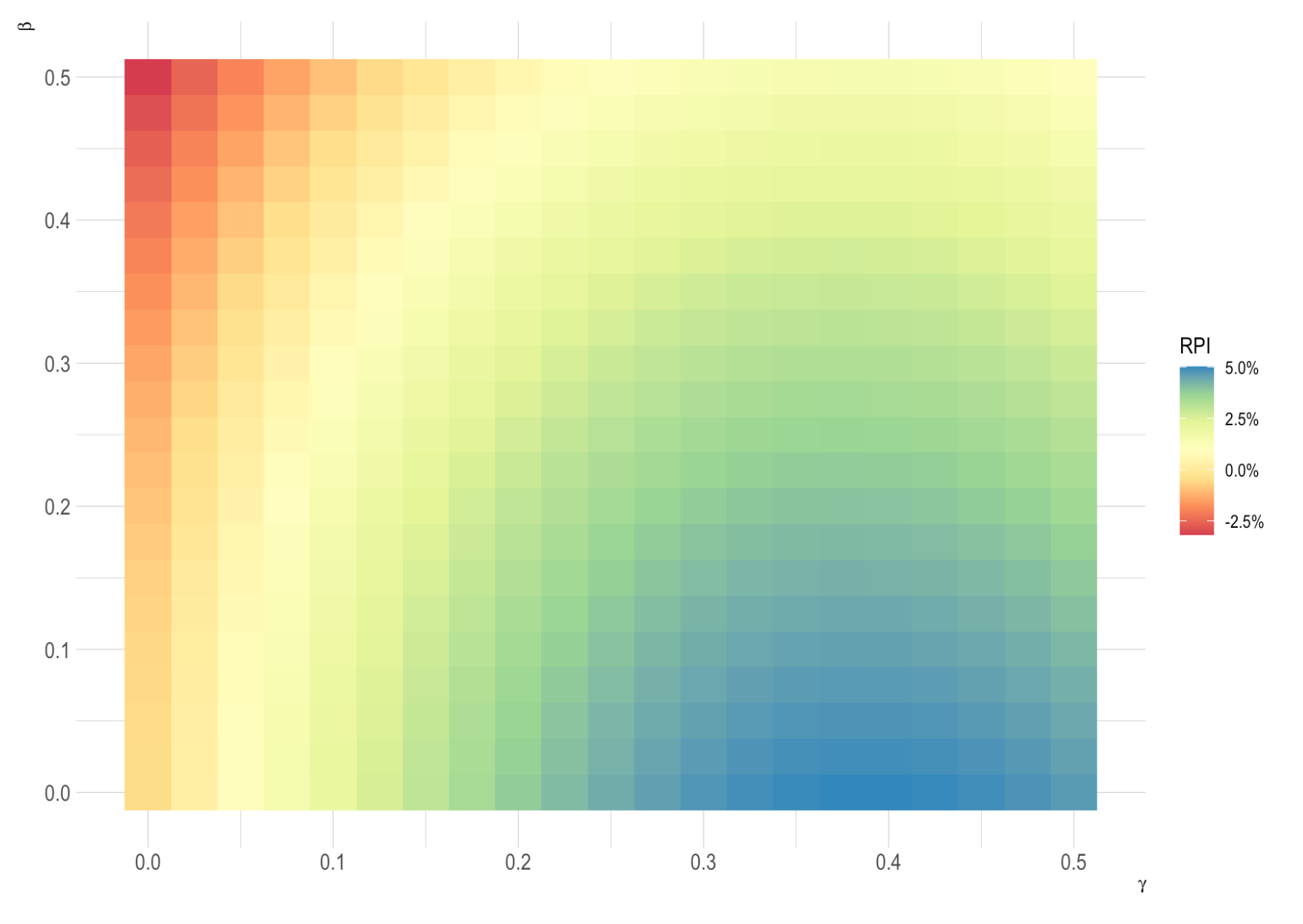}
\end{figure}

\begin{figure}[h!]
 \centering
 \caption{RPI heatmap for the over-parametrised case
 ($\tau = 0.7$, $t=200$)}
\includegraphics[width=0.8\textwidth]{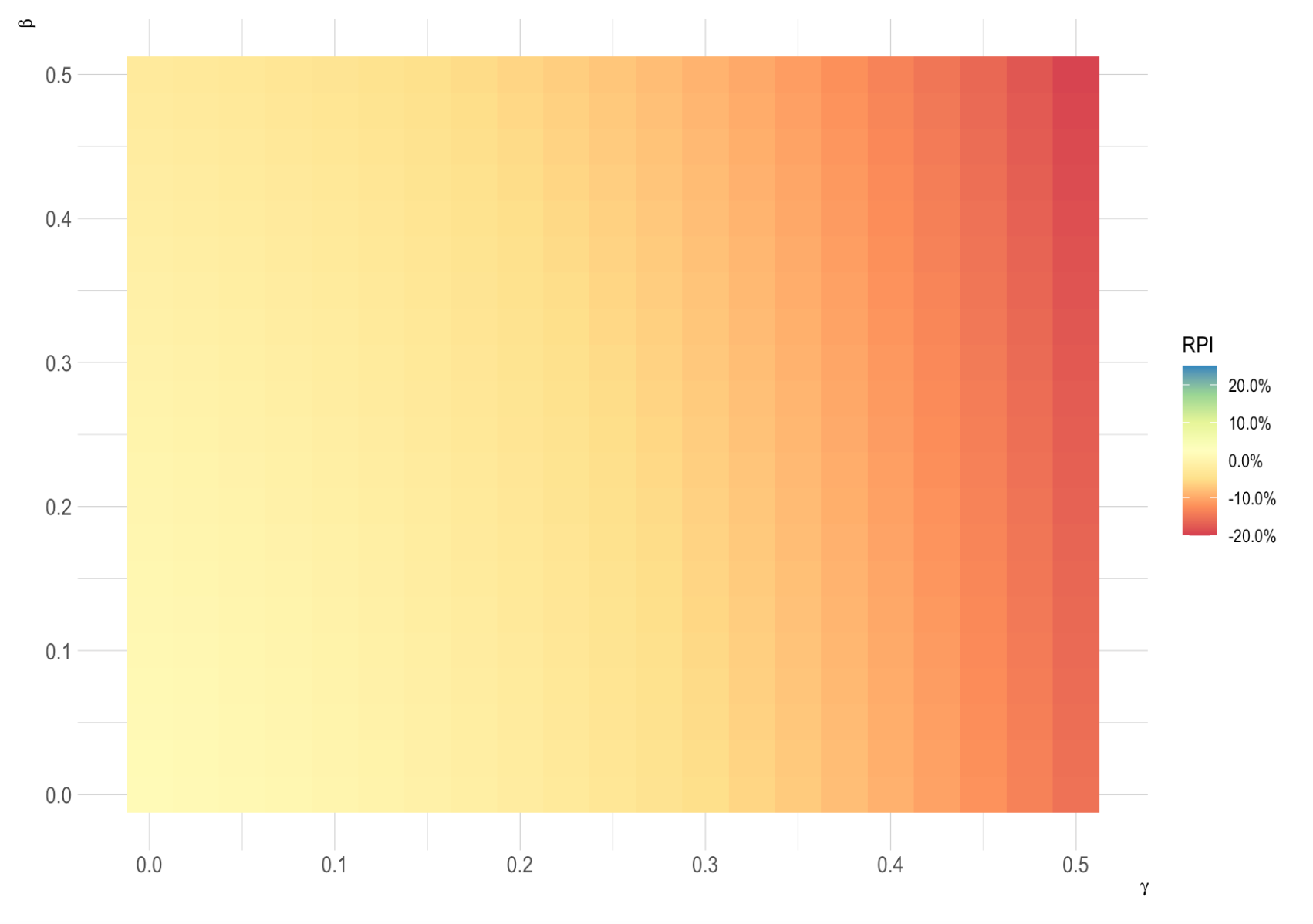}
\end{figure}
\end{document}